\documentclass[12pt]{article}
\usepackage{amsmath, amssymb, amsthm}
\usepackage{enumerate}
\newtheorem{theorem}{Theorem}[section]
\newtheorem{lemma}{Lemma}[section]
\newtheorem{proposition}{Proposition}[section]

\newtheorem{corollary}{Corollary}[proposition]
\newtheorem*{theorema}{Theorem A}
\newtheorem*{theoremb}{Theorem B}
\newcommand{\taf}{{\hskip 5pt} $\blacksquare$
                  \renewcommand{\qedsymbol}{}}
\begin{document}
\title{Spaces between $H^{1}$ and $L^{1}$}
\author{Wael Abu-Shammala and Alberto Torchinsky}
\date{}
\maketitle
\begin{abstract} In this paper we consider the spaces $X_s$
 that lie between $H^1(R^n)$ and $L^1(R^n)$. We discuss their
 interpolation
properties  and  the
 behavior of maximal functions and
singular integrals acting on them.
\end{abstract}
Because of their similarities, but mainly because of their
differences, it is a matter of interest to determine the
relationship between the Hardy space $H^1(R^n)$ and the space of
integrable functions $L^1(R^n)$. The purpose of this paper is to
address some unanswered questions concerning the family of spaces
$X_s(R^n)$ that lie between $H^1$ and $L^1$, and thus gain a
better understanding of the gap that separates them.

The spaces $X_s$ were  introduced by Sweezy, see \cite{sweezy}.
They form a nested family that starts
 at $H^1= X_1$ and approaches $L^1_0$, the subspace of $L^1$
functions with vanishing integral, as $s\to\infty$.  Here we
consider the whole range of $X_s$ spaces. First, $X_s=H^1$ for
$0<s\le 1$; also, $X_{\infty}=L^{1}_0$, see \cite{AT}.  Further,
we  show that, for $f\in X_s$,
\[K(t,f;H^1,L^1)\le c\,\min(t,t^{1/s'})\,\|f\|_{X_s}\,.\]
This estimate gives that $X_s$ is continuously embedded in the
Hardy-Lorentz space $H^{1,r}$  consisting of those distributions
with   non-tangential maximal function in the Lorentz space
$L^{1,r}$, for $1<s<r\le\infty$. Therefore, also the spaces
$H^{1,s}\cap L^1$, $1\le s< \infty$, form a nested family of
subspaces of $L^1$ that increase from $H^1$ to $L^1$. As for
Calder\'{o}n-Zygmund singular integral operators, they   map $X_s$
into $L^{1,r}$ for $1<s<r\le \infty$.
 We conclude the paper by introducing the closely related family
of $X^s$ spaces, $0< s\le\infty$, that increases towards $L^1$, and
then showing how $X_s$ and $X^s$ atoms can be used to build other
spaces, including  analogues of the local spaces considered
  in \cite{LY}, that  lie between $H^1$ and $L^1$.
\section{Atomic decompositions in Banach spaces}
Our first result is of general nature and will ensure that the
various atomic spaces considered below are indeed Banach spaces.

Let $\mathcal{A}$
 be a non-empty subset of the unit ball of a Banach space $(B,\|\cdot\|)$.
 The atomic space $A$
 spanned by $\mathcal{A}$ consists of all those $\varphi\in B$
   of the form
\[\varphi=\sum \lambda_ja_j\,,\quad \sum|\lambda_j|<\infty\,,a_j\in
\mathcal{A}\,.\] It is readily seen that, endowed with the atomic
norm $\|\varphi\|_A=\inf\big\{\sum_1^{\infty}|\lambda_j|:
\varphi=\sum_1^{\infty}\lambda_j\,a_{j} \,\big\}$,  $A$ becomes a
normed space. In fact, it is also complete.
\begin{lemma}
The atomic space $(A, \|\cdot\|_A)$ is a Banach space.
\end{lemma}
\begin{proof} Since for
$\varphi=\sum_1^{\infty}\lambda_j\,a_{j}$   we have
$\|\varphi\|\le\sum_1^{\infty}|\lambda_j|$, it readily follows that
$\|\varphi\|\le \|\varphi\|_A$, and %$A\hookrightarrow B$.
$A$ is continuously embedded in $B$.

To verify that
 $(A, \|\cdot\|_A)$ is complete, it suffices to
show that if $\{\varphi_n\}$ is a
 sequence of elements in $A$  such that
$\sum_1^{\infty}\|\varphi_n\|_A<\infty$, then the sum converges in
$A$, i.e., for some $\varphi\in A$, $\lim_{N\to\infty}\big\|\varphi
-\sum_{n=1}^N\varphi_n\big\|_A=0\,$. First observe that, since also
$\sum_1^{\infty}\|\varphi_n\|<\infty$,  the sum
 converges to some $\varphi$ in $B$. We will show that $\varphi\in A$,
 and that the sum also converges to $\varphi$ in $A$.

Let $\varphi_n=\sum_{j=1}^{\infty}\lambda_{j,n}\,a_{j,n}$, where the
$\lambda_{j,n}$'s satisfy $\sum_{j=1}^{\infty}|\lambda_{j,n}|\le \,2
\|\varphi_n\|_A$, $n=1,2,\cdots$ Having fixed these decompositions,
we may restrict our attention to the countable set of atoms
$\{a_{j,n}\}$. So, we rename these elements  $\{a_j\}$ and, by
adding zeroes to the original $\lambda_{j,n}$'s as needed, we have
$\varphi_n=\sum_{j=1}^{\infty}\lambda_{j,n}a_j$, all $n$. Clearly
$\sum_{n=1}^{\infty}\sum_{j=1}^{\infty}|\lambda_{j,n}|<\infty$.
Moreover, if $\mu_j=\sum_{n=1}^{\infty} \lambda_{j,n}$, then
$\sum_j|\mu_j|\le
\sum_{j=1}^{\infty}\sum_{n=1}^{\infty}|\lambda_{j,n}|<\infty$. Now,
since $\varphi=\sum\varphi_n=\sum_{j=1}^{\infty}\mu_j\,a_j$,
 $\varphi\in A$. Finally, given $\varepsilon>0$,
let $N_0$ be such that
$\sum_{n=N}^{\infty}\sum_{j=1}^{\infty}|\lambda_{j,n}|\le\varepsilon$,
for $N\ge N_0$.  Then,  for $N\ge N_0$,
\begin{align*}
\bigg\| \varphi-\sum_{j=1}^{\infty}\bigg(\sum_{n=1}^{N-1}
\lambda_{j,n}\bigg)a_j\bigg\|_A &=\bigg\|\sum_{j=1}^{\infty}
\bigg(\sum_{n=N}^{\infty}\lambda_{j,n}\bigg)a_j\bigg\|_A\\
&\le \sum_{n=N}^{\infty}\sum_{j=1}^{\infty} |\lambda_{j,n}|\le\,
\varepsilon\,,
\end{align*}
and we have finished. \taf
\end{proof}
\section{The spaces $X_s$}
 For an exponent $1<q\le 2$ with conjugate
 $2\le p<\infty$,
$1/p+1/q=1$, and $0<s\le \infty$, we say that a compactly supported
function $a$ with vanishing integral is a $(q,s)$ atom with defining
cube  $Q$ if
\[{\text{supp$(a)\subseteq Q$}},\quad
\int_Q a(x)\,dx=0\,,\quad p^{1/s}|Q|\left(\frac1{|Q|}
\int_Q|\,a(x)\,|^qdx\right)^{1/q}\le 1\,.\] When $s=\infty$,  $a$ is
a usual $L^q$ $1$-atom.

We denote by ${\mathcal A}_s$ the collection of $(q,s)$ atoms,
$0<s\le\infty$.  Now, since for $a\in {\mathcal A}_s$ we have
\[\int_Q|a|\le
|Q|\left(\frac1{|Q|}\int_Q|a(x)|^q\,dx\right)^{1/q} \le
p^{-1/s}\le 1\,,\] ${\mathcal A}_s$ is contained in the unit ball
of $L^1$, and Lemma 1.1 applies. The resulting Banach space is
$X_s(R^n)=X_s$, the space introduced by Sweezy, who also showed
that $X_1=H^1$, see \cite {sweezy}.

As a first step in determining the  relationship between the
$X_s$'s for the different values of $s$ we have
\begin{lemma} Suppose $a$ is a $(q,s)$ atom. Then,
\[\|a\|_{X_u} \le p^{1/u-1/s}\,,\quad 0<u\le\infty\,.\]
\end{lemma}
\begin{proof} The conclusion
follows readily since, for a $(q,s)$ atom $a$,
\[p^{1/u}|Q|\left(\frac1{|Q|}\int_Q \bigg[\,\frac{|a(x)|}
{p^{1/u-1/s}}\bigg]^q dx\right)^{1/q}\le 1\,\,.{\hskip
10pt}\blacksquare\]
\renewcommand{\qedsymbol}{}
\end{proof}
From Lemma 2.1 we get $\|f\|_{X_s}\le \|f\|_{X_r}$,   $0< r\le
s\le\infty$, and the $X_s$'s  are nested. Also,  the norm in $X_s$
reduces to the atomic $H^1$ norm for $0<s<1$.
\begin{proposition} Suppose $0<s<1$. Then $X_s=H^1$, with equivalent
norms.
\end{proposition}
\begin{proof} We have already noted that $\|f\|_{H^1}\le
\|f\|_{X_s}$.
  Let now $f\in H^1$ have an atomic
decomposition $f=\sum_j\lambda_j\,a_j$ in terms of
 $L^{\infty}\, 1$-atoms  $a_j$. Then,
\[2^{1/s}|Q|\left(\frac1{|Q|}\int_Q\bigg[\,\frac{|a_j(x)|}
{2^{1/s}}\bigg]^2 dx\right)^{1/2}\le 1\,,\] $a_j(x)/ 2^{1/s}$ is a
$(2,s)$ atom, and $\|f\|_{X_s}\le 2^{1/s}\sum_j |\lambda_j|\,.$
Taking the infimum over all possible decompositions of $f$ it
follows that $\|f\|_{X_s}\le 2^{1/s}\|f\|_{H^1}$, and we have
finished.\taf
\end{proof}
The situation is different for $s>1$. As the $q_j$'s approach $1$,
 the $p_j$'s tend to $\infty$,  the sums escape $H^1$, and
$X_s\ne H^1$. In fact, $X_s$ contains strictly $X_r$ for $1\le
r<s$,  and, although $X_r$ is densely embedded in $X_s$ for $r<s$,
it is of first category in $(X_s,\|\cdot\|_{X_s})$.

Concerning  $s$ large, let $X=\bigcup_{s<\infty}\!X_s$. We
%$\varinjlim X_{s}$. We
introduce a topology in $X$ that is easier to deal with than the
inductive topology there. For $f\in X$, let
 $\|f\|_X=\lim_{s\to\infty} \|f\|_{X_s}$.
It is not hard to see that $\|\cdot\|_X$ is a norm. The homogeneity
and triangle inequality follow easily. Moreover, if $ \|f\|_{X}=0$,
 $ \|f\|_{1}=0$, and $f=0$ a.e.  In fact, for each $s$, the
inclusion mapping is continuous from $X_s$ to $(X, \|\cdot\|_X)$,
and, consequently,   $\varinjlim X_{s}$ is also continuously
included in $X$. Similarly,  $(X, \|\cdot\|_X)$ is continuously
embedded in $(X_{\infty},\|\cdot\|_{X_{\infty}})$, but  $X$ and
$X_{\infty}$ are not the same space. To see this  note that $X_r$
is of first category in $X_{\infty}$ for $r<\infty$, and the same
is true for $X$: if ${\cal U}$ is  open  in $X_{\infty}$,  it can
not be open in any $X_r$, and hence it is not open in $X$.
Finally, $X_{\infty}=L^1_0$, see \cite{AT}.
\subsection{$X_s$ as an intermediate space between $X_{s_1}$ and
$X_{s_2}$} Recall that  the $K$ functional of $f\in
X_{s_1}+X_{s_2}$ at $t>0$ is defined by
\[K(t,f;X_{s_1},X_{s_2})=\inf_{f=f_1+f_2}\|f_1\|_{X_{s_1}}+t\,\|f_2\|_{X_{s_2}}\,,\]
where $f=f_1+f_2$, $f_1\in X_{s_1}$ and $f_2\in X_{s_2}$. We begin
by estimating the $K$ functional for $f\in X_s$, $1\le s_1< s<
s_2\le \infty$. The reader
  will have no difficulty in verifying that a
similar result holds with $X$ in place of $X_{\infty}$.
\begin{lemma} Given $1\le s_1<s<s_2\le \infty$, let $0<\eta<1$ be given
by $1/s=(1-\eta)/s_1+\eta/s_2$.
 Then, for $f\in X_s$,
\[K(t,f;X_{s_1},X_{s_2})\le \min\,(\,t\,,\,t^{\eta})\,\|f\|_{X_s}.\]
\end{lemma}
\begin{proof}
 Since $X_{s} \hookrightarrow X_{s_2}$,
 $K(t,f;X_{s_1},X_{s_2})\le t \,\|f\|_{X_{s_2}}\le t\,\|f\|_{X_s}$.
 This estimate suffices for $t$  small.

Suppose now that $t$ is large, $t>1$, say, and let $\alpha>0$ be
given by $1/\alpha=1/s_1 -1/s_2$. Let $f\in X_s$ have an atomic
decomposition $f=\sum_j \lambda_j a_j$ in terms of  $(q_j,s)$
atoms $a_j$, and let $p_j$ denote the conjugate exponent to $q_j$.
Finally, put ${\mathcal J}_1=\{j: p_j\le t^{\alpha}\}$ and
${\mathcal J}_2=\{j: p_j > t^{\alpha}\}$. By Lemma 1.2 we have
\begin{align*}\bigg\|\sum_{j\in {\mathcal J}_1}\lambda_j\,a_j\,\bigg\|_{X_{s_1}}
&\le \sum_{j\in {\mathcal J}_1}| \lambda_j| \,
\|a_j\|_{X_{s_1}}\\
&\le \sum_{j\in {\mathcal J}_1}|\lambda_j|\, p_j^{1/s_1-1/s}\le
t^{\alpha(1/s_1-1/s)}\, \sum_{j\in {\mathcal J}_1}|\lambda_j|\,,
\end{align*}
\begin{align*}\bigg\|\sum_{j\in {\mathcal J}_2}\lambda_j\,a_j\,\bigg\|_{X_{s_2}}
&\le \sum_{j\in {\mathcal J}_2}| \lambda_j| \,
\|a_j\|_{X_{s_2}}\\
&\le \sum_{j\in {\mathcal J}_s}|\lambda_j|\, p_j^{1/s_2-1/s}\le
t^{\alpha(1/s_2-1/s)}\, \sum_{j\in {\mathcal J}_2}|\lambda_j|\,.
\end{align*}
Now, since, as is readily seen, $\alpha(1/{s_1}-1{s})=\eta$ and
$\alpha(1/s_2-1/s)=\eta-1$, we get
\[K(t,f;X_{s_1},X_{s_2})\le \bigg\|\sum_{j\in {\mathcal
J}_1}\lambda_j\,a_j\,\bigg\|_{X_{s_1}}+t\,\bigg\|\sum_{j\in
{\mathcal J}_2}\lambda_j\,a_j\,\bigg\|_{X_{s_2}}\le
t^{\eta}\,\sum_{j}|\lambda_j|\,.\] Thus, taking the infimum over
the decompositions of $f$ in $X_s$, it follows that
\[K(t,f;X_{s_1},X_{s_2})\le t^{\eta}\,\|f\|_{X_s}\,.\]
 The
conclusion now obtains by combining the estimates for $t$ small and
$t$ large.\taf
\end{proof}
\begin{corollary}Let $1<s<\infty$.
 Then, for $f\in X_s$,
\[K(t,f;H^1,L^1)\le c\,\min\,(\,t\,,\,t^{1/s'})\,\|f\|_{X_s}.\]
\end{corollary}
\begin{proposition}
Let $f\in X_s$, $1<s<\infty$. Then $f$ is in the Hardy-Lorentz
space $H^{1,r}$, $r>s$, and $\|f\|_{H^{1,r}}\le c_r\,
\|f\|_{X_s}$, $c_r= O\big(c/(1/s-1/r)\big)$.
\end{proposition}
\begin{proof}  Let $f\in X_s$. We will show that the non-tangential
 maximal function $Nf$ is in  $L^{1,r}$ for
 $r>s$. Since the non-tangential maximal function  of a function
in $H^{1}$ is   in $L^{1}$, and  that of a function in $L^{1}$ is
in $L^{1,\infty}$,  by elementary interpolation considerations and
Corollary 2.1.1 it follows that
\[K(t,Nf;L^1,L^{1,\infty})\le c\, K(t,f;H^1,L^1)\le c\,
\min(t,t^{1/s'})\,\|f\|_{X_s}.\] Given $1<s<r<\infty$, let
$1/s'<\theta<1$ be chosen so that $1/r={1-\theta}$.
 Integrating the above inequality we get
\begin{align*}\|Nf\|^s_{(L^{1},L^{1,\infty})_{\theta,s}} &\le c\,
\left(\int_0^{\infty}\left(\frac{\min\,(t,
\,t^{1/s'})}{t^{\theta}}\right)^s\frac{dt}{t}\right)\,{\|f\|^s_{X_s}}\\
&\le c_{r}\,{\|f\|^s_{X_s}}.
\end{align*}
Clearly, $c_{r}\le c/(1/s-1/r)$.
 Furthermore, since $s<r$, we also have, see \cite{BL},
\[\|Nf\|_{1,r}\sim
\|Nf\|_{(L^{1},L^{1,\infty})_{\theta,r}}\le c\,
\|Nf\|_{(L^{1},L^{1,\infty})_{\theta,s}},\] and the conclusion
follows by combining the two estimates.\taf
\end{proof}
 Lemma  2.2   applies   to Calder\'{o}n-Zygmund singular integrals, and other
  operators, such as  the Marcinkiewicz integral, see \cite{DLX},
  that map $H^1$ into  $L^1$,
 and $L^1$ into weak $L^1$.  Thus, these operators  also map $X_s$
 into the Lorentz space $L^{1,r}$, for $1<s<r\le \infty$. It also
 applies to some
Calder\'{o}n-Zygmund singular integral operators with rough kernels
that are known to be of weak-type $(1,1)$, see \cite{S}, and to map
$H^1$ into $L^{1,2}$, see \cite{ST}.  Lemma 2.2 then gives that they
also map $X_s$ into $L^{1,r}$ for $r>2s$.
\section{Concluding remarks}
In order to reach $L^1$ from $H^1$, one more atom, this one with
nonvanishing integral, needs to be added to the families ${\mathcal
A}_s$; the characteristic function of $Q_1$,
%$=[\,-1/2,1/2\,]\times\cdots\times [\,-1/2,1/2\,]$,
the cube  of sidelength $1$ centered at the origin, will do. Let
$X^s$ denote the Banach space spanned by ${\mathcal A}_s\cup
\chi_{Q_1}$, $1\le s\le\infty$, and note that if $f\in X^s$ has an
atomic decomposition of the form $f=\sum_j\lambda_j\,a_j + \lambda
\,\chi_{Q_1}$, $\lambda$ is uniquely determined and is equal to
$\int_{R^n} f$. Thus, for  $f$ in $X^s$ we have $\|f\|_{X^s}= \inf
\sum_j|\lambda_j|+|\lambda|$, where the infimum is taken over the
atomic decompositions
 of $f$. In other words,  $X^{s}= X_s+\, {\rm sp}(\chi_{Q_1})$, in the sense of
sum of Banach spaces. Moreover, the family $X^{s}$ is nested and
reaches $L^1$. Other than this important property, the family
$X^{s}$ behaves very much like $X_{s}$, and, consequently, we only
give the description of its dual. It is the space BMO$^{s}$
consisting of those functions $\varphi(x)$ such that
\[
A(\varphi)=\sup_{p>1}\frac{1}{p^{1/s}}\,\sup_{Q}\left(\frac{1}{ |
Q|}\int\nolimits_{Q}| \varphi (x)-\varphi _{Q}|
^{p}dx\right)^{1/p}<\infty\,,
\]
\[
B(\varphi )=\Big|\int_{Q_1}\varphi (x)\,dx\Big|<\infty\,,
\]
normed with $\| \varphi \|_{BMO^{s}}=\max (A(\varphi ),B(\varphi
))$, $1\le s<\infty$. Of course, the dual of $X^{\infty}$ is
$L^{\infty}$, see \cite{AT}.

Finally, one can define spaces analogous to the local version of
$H^1$ spaces at the origin introduced in \cite{LY}, with the $H^1$
atoms there replaced by $X_s$ atoms.   Let $Q_{\delta}$
 denote the cube   of sidelength $\delta$
centered at the origin. The family ${\mathcal C}_s$ of central
$(q,s)$ atoms consists of those  $a\in {\mathcal A}_s$ with defining
cube $Q_{\delta}$ for some $\delta>0$. The atomic space generated by
${\mathcal C}_s$ is denoted  $HX_s$, for $1\le s \le \infty$.

$HX_1$ is a dense subset of $H^1$ which is embedded continuously
in $H^1$ and, as we will see below,  it is not $H^1$. The
 $HX_s$'s form a nested family of subspaces of $L^1_0$ and, for
each $s$,   $HX_s$ is a dense subset of $X_s$  continuously
embedded in $X_s$.

The dual of $HX_s$ is the space CMO$(s)$ which consists of those
functions $\varphi$ such that
\[\|\varphi\|_{CMO(s)}=\sup_{p>1}\frac{1}{p^{1/s}}\,\sup_{\delta>0}
\left(\frac1{|Q_{\delta}|}\int_{Q_{\delta}}|\varphi(x)-
\varphi_{Q_{\delta}}|^pdx\right)^{1/p}<\infty\,,\quad s<\infty\,,
\]
\[\|\varphi\|_{CMO(\infty)}=\sup_{p>1}\,\sup_{\delta>0}
\left(\frac1{|Q_{\delta}|}\int_{Q_{\delta}}|\varphi(x)-
\varphi_{Q_{\delta}}|^pdx\right)^{1/p}<\infty\,,\quad s=\infty\,.
\]
Now,  CMO$(1)$ strictly contains BMO. To see this consider, for
$n=1$,
\[
\varphi(x) = \begin{cases} \displaystyle  \ln |x+2|,
& -\infty<x<-2\\%[5pt]
\displaystyle 0, & -2\le x\le -1\\
\displaystyle   \ln |x+1|, & -1<x<\infty\,.
\end{cases}
\]
$\varphi$ is in CMO$(1)$ but not in BMO, which shows that, unlike
$X_1$, $HX_1$ is not $H^1$. The space CMO$(\infty)$ coincides with
$L^{\infty,*}$, the dual of $L^1_0$, as a set of functions, but
with a weaker norm.

Finally,  let $HX^s$ denote the Banach space  spanned by ${\mathcal
C}_s \cup \chi_{Q_1}$.  The  $HX^s$'s form an increasing family of
subspaces of $L^1$ and $HX^s$ contains strictly $HX_s$.
$HX^{\infty}$ is close to $L^1$, so the $HX^s$'s cover the gap left
by the $HX_s$'s. The dual of $HX^s$ is CM${\rm O}^s$, which consists
of all functions $\varphi$ in $CMO(s)$, but normed with
$\|\varphi\|_{\rm{CMO}^s}=\|\varphi\|_{{\rm {CMO}}(s)}+\big|
\int_{Q_1}\varphi(x)\,dx\big|$, $1\le s<\infty$.

DEPARTMENT OF MATHEMATICS, INDIANA UNIVERSITY,
\\BLOOMINGTON, IN 47405
\\{\it E-mail:} wabusham@indiana.edu, torchins@indiana.edu

\begin{thebibliography}{99}
%\bibitem{AST} {W. Abu-Shammala, J.-L. Shiu,} and {A.
%Torchinsky}, {\em Characterizations of the Hardy space $H^1$ and
%{\rm BMO}}, preprint.
\bibitem{AT} {W. Abu-Shammala} and {A.
Torchinsky}, {\em The atomic decomposition in $L^1(R^n)$}, Proc.
Amer. Math. Soc., to appear.

\bibitem{BL} {J. Bergh} and {J. L${\ddot{\rm o}}$fstr${\ddot{\rm
o}}$m}, {\em Interpolation spaces, an introduction},
Springer-Verlag, 1976.

\bibitem{DLX} {Y. Ding, S. Z. Lu,} and {Q. Xue}, {\em Marcinkiewicz integral
on Hardy spaces,} Integr. Equ. Oper. Theory {\bf 42}, (2002),
174-182.

\bibitem{Fefferman_Riviere} {C. Fefferman, N. M. Rivi${\grave{\rm
e}}$re}, and {Y. Sagher}, {\em Interpolation between $H^p$ spaces:
the real method}, Trans.   Amer. Math. Soc. {\bf 191} (1974),
75--81.

%\bibitem{fridli} {S. Fridli}, {\em Transition from the dyadic to the
%real nonperiodic Hardy space},  Acta Math. Acad. Paedagog.
%Nih\'{a}zi (N.S.) {\bf 16} (2000), 1--8, (electronic).

\bibitem{GCRdF} {J. Garc\'{i}a-Cuerva} and {J. L. Rubio de
Francia}, {\em Weighted norm inequalities and related topics}, Notas
de Matem\'{a}tica {\bf 116}, North Holland, 1985.

 \bibitem{LY} {S. Z. Lu} and {D. C. Yang}, {\em The local version
of $H^p(R^n)$ spaces at the origin},
 Studia Math.  {\bf 116} (1995), 103--131.

\bibitem{S} {A. Seeger}, {\em Singular integral operators with rough
convolution kernels,} J. Amer. Math. Soc. {\bf 9} (1966), 95-105.

 \bibitem{ST} {A. Seeger} and {T. Tao}, {\em Sharp Lorentz space
 estimates for rough operators}, Math. Ann. {\bf 320} (2001),
 381-415.

\bibitem{sweezy} {C. Sweezy}, {\em Subspaces of  $L^1(R^d)$},
Proc.   Amer. Math. Soc. {\bf 132} (2005), 3599--3606.

\bibitem{torchinsky} {A. Torchinsky},
{\em Real-variable methods in harmonic analysis},  Dover
Publications, Inc., 2004.
\end{thebibliography}
\end{document}